  \documentclass[]{amsart}
   \newtheorem{lemma}{Lemma}[section]
   \newtheorem{theorem}[lemma]{Theorem}

   \newtheorem{coro}[lemma]{Corollary}
   \newtheorem{definition}[lemma]{Definition}
   \newcommand{\eps}{\varepsilon}
   \newcommand{\Hsob}{\smash{{\stackrel{\circ}{H}}}^1(D)}

\newcommand{\om}{\omega}
\newcommand{\Om}{\Omega}
\newcommand{\D}{\Delta}
\newcommand{\p}{\partial}
\renewcommand{\phi}{\varphi}

\title[Ergodic Dynamics of Geophysical Flows]
{Ergodic Dynamics of the Coupled Quasigeostrophic Flow-Energy
Balance System}

\author{Aijun Du, Jinqiao Duan}
\address[A. Du and J. Duan]
{Department of Applied Mathematics\\
 Illinois Institute of Technology\\
   Chicago, IL 60616, USA }
\email[]{duaijun@iit.edu; duan@iit.edu}

\author{Hongjun Gao}
\address[H. Gao]
{Department of Mathematics\\
Nanjing Normal University\\
Nanjing 210097, China} \email[H.~Gao]{gaohj@njnu.edu.cn}

\author{Tamay \"{O}zg\"{o}kmen }
\address[T. \"{O}zg\"{o}kmen]
{RSMAS/MPO\\University of Miami\\Miami, Florida, USA} \email[T.
 \"{O}zg\"{o}kmen]{tamay@rsmas.miami.edu}

\date{May 20, 2004}


\begin{document}

\begin{abstract}

 The authors consider a mathematical model for the coupled
atmosphere-ocean system, namely, the coupled quasigeostrophic
flow-energy balance model. This model consists of the large scale
quasigeostrophic oceanic flow model and the transport equation for
oceanic temperature, coupled with an atmospheric energy balance
model. After   reformulating this coupled model as a random
dynamical system (cocycle property), it is shown that the coupled
quasigeostrophic-energy balance fluid system has a random
attractor, and under further conditions   on the physical data and
the covariance of the noise, the system is ergodic, namely, for
any observable of the coupled atmosphere-ocean flows, its time
average approximates the statistical ensemble average, as long as
the time interval is sufficiently long.

\end{abstract}

\maketitle

\section{Mathematical model}\label{sect1}

We consider large scale geophysical flows modeled by the
quasigeostrophic flow equation in the horizontal $xy-$plane,
  in terms of vorticity
$q(x,y,t)$, and the transport equation for the oceanic temperature
$T(x,y, t)$, coupled with the atmospheric energy balance equation
proposed by North and Cahalan \cite{NorCah81}  for the air
temperature $\Theta(x,y,t)$, on the domain $D=\{ (x,y): 0\leq x
\leq l, 0 \leq y \leq l \}$:

\begin{equation}\label{q1}
\begin{split}
\Theta_t   =&  \Delta\Theta -(a +\Theta)+S_a(x,y) - b(y)(S_o(x,y)+\Theta -T(x,y))  + \dot{w},  \\
\frac{\partial q}{\partial t}
          = &\nu \D  q  - r q + Pr\;Ra\; \partial_y T - J(\psi, q + \beta  y ),\\
\frac{\partial T}{\partial t}  =& \Delta T - J(T, \psi),
\end{split}
\end{equation}

where $\psi(x,y,t)$ is the stream function, $\beta \geq 0$ is the
meridional gradient of the Coriolis parameter, $\nu >0$ is the
viscous dissipation constant, and $r>0$ is the Ekman dissipation
constant. Furthermore,
\[
q(x,y, t)= \D \psi(x,y,t)
\]
is the vorticity,  $a>0$ is a   constant parameterizing the effect
of the earth's longwave radiative cooling, $b(y)$ is the
latitudinal fraction of the earth covered by the ocean basin, Pr
is the Prandtl number, and Ra is the Rayleigh number. Note that
$S_a(x,y)$ and $S_o(x,y)$ are empirical functions representing the
effects (on atmosphere and ocean, respectively) of the shortwave
solar radiation. Moreover, $J(g,h)=g_xh_y-g_yh_x$ is the Jacobian
operator and $\D=\p_{xx}+\p_{yy}$ is the Laplacian operator. All
these equations are in non-dimensionalized forms. The fluctuating
noise $\dot{w}(x,y, t)$   is usually  of a shorter time scale than
the response time scale of the   air  mean temperature. So we
neglect the autocorrelation time of this fluctuating process and
thus assume that the noise is white in time. The spatially
correlated white-in-time noise $\dot{w}(x,y, t)$ is described as
the generalized time derivative of a Wiener process $w(x,y, t)$
defined in a probability space $(\Omega, \mathbb{F}, \mathbb{P})$,
with mean zero and covariance operator $Q$. For geophysical
background of similar coupled atmosphere-ocean models, see
\cite{Dij00,CheGhi96}. Recently, a few authors have considered the
randomly forced quasigeostrophic equation, in order to incorporate
the impact of uncertain geophysical forces (
\cite{Griffa},
\cite{Holloway},
\cite{DelSole-Farrell}).

The  fluid boundary condition is no normal flow and free-slip on
the whole boundary:
\begin{equation}
\label{bvorticity} \psi=0,  \;  q = 0.
\end{equation}

The flux boundary conditions are assumed for the ocean temperature
$T$ and air temperature $\Theta$:
\begin{equation}
\label{temp}
   \frac{\partial\Theta}{\partial n}  = \frac{\partial T }{\partial n} = 0,
 \end{equation}
with $n$ is the outer unit normal vector on boundary. By equation
(\ref{q1}) and the boundary condition, we find that $\int_D T(x,
y) dxdy = \mbox{constant}$, and thus we assume that $\int_D T(x,
y) dxdy = 0$ without loss of generality.

This is a coupled system of both deterministic and stochastic
partial differential equations.

This paper is organized as follows.
 In the next section, we
introduce basic concepts in random dynamical systems. Then we
reformulate the coupled fluid system (\ref{q1}) as a random
dynamical system with help of a Ornstein-Uhlenbeck random
stationary process in \S 3. After obtaining basic estimates for
the system in \S 4, we show in \S 5 that the coupled fluid system
(\ref{q1}) admits a random attractor. Under further conditions on
the physical data and the covariance of the noise, we show that
the system is actually ergodic, namely, for any observable of the
coupled atmosphere-ocean flows, its time average approximates the
statistical ensemble average, as long as the time interval is
sufficiently long.

\section{Random dynamical systems}\label{s2}

In order to investigate the long time dynamics of  the coupled
fluid system (\ref{q1})
under the influence of random forces, we need some appropriate concepts and tools from the theory of {\em random dynamical systems}.\\

A random dynamical system   consists of two components. The first
component is a {\em driven   flow} $(\Omega,{\mathcal
F},\Bbb{P},\theta)$ as a model for   noise, where
$(\Omega,{\mathcal F},\Bbb{P})$ is a probability space and
$\theta$ is a $\mathcal{F}\otimes {\mathcal B}(\Bbb{R}),{\mathcal
F}$ measurable flow: we have
\[
\theta_0={\rm id},\qquad \theta_{t+\tau}=\theta_t\circ\theta_\tau
=:\theta_t\theta_\tau
\]
for $t,\,\tau\in\Bbb{R}$.  To express that the noise is stationary
and ``chaotic",  the measure $\Bbb{P}$ is supposed to be ergodic
with respect to $\theta$. The second  component of a random
dynamical system is a ${\mathcal B}({\Bbb R}^+) \otimes{\mathcal
F}\otimes{\mathcal B}(H),{\mathcal B}(H)$-measurable mapping
$\phi$ satisfying the {\em cocycle} property
\[
\phi(t+\tau,\omega,x)=\phi(t,\theta_\tau\omega,\phi(\tau,\omega,x)),\qquad
\phi(0,\omega,x)=x,
\]
where the phase space $H$ is a separable metric space and $x$ is
chosen arbitrarily in $H$. We will denote this
random dynamical system by symbol $\phi$.\\
A standard model for   a spatially correlated white-in-time noise
 is the generalized time derivative of a two-sided {\em Brownian motion or Wiener process}
 $w(x, y, t)$. Let $U$ be a
separable Hilbert space with  scalar product $(\cdot, \cdot)$ and
 the induced norm $\| \cdot \|$. For stochastic partial
differential equations containing such a noise, an appropriate or
canonical probability space is
\[
(C_0({\Bbb R},U),{\mathcal B}(C_0({\Bbb R},U)),\Bbb{P}),
\]
where $C_0(\Bbb{R},U)$ is the   space of continuous functions on
$\Bbb{R}$, which take zero value at time zero. This space is given
the compact-open topology (i.e., uniform convergence on compact
intervals in $\Bbb{R}$). This topology is metrizible as it can be
generated by the complete metric
 $$
 d(g_1, g_2) = \sum_{n=1}^{\infty} \frac1{2^n}
 \frac{d_n(g_1, g_2)}{1+d_n(g_1, g_2)},
 $$
 where $d_n(g_1, g_2)=\max_{|t|\leq n} \|g_1(t)-g_2(t)\|$, for $g_1,
 g_2$ in $C_0({\Bbb R},U)$. Thus we have open balls or open sets
 in $C_0({\Bbb R},U)$, and
  ${\mathcal B}(C_0({\Bbb R},U))$ is the corresponding
Borel $\sigma$-algebra. Suppose   the Wiener process $w$ has
covariance operator $Q$ on $U$. Let $\Bbb{P}$ denote  the {\em
Wiener measure} with respect to $Q$. Note that $\Bbb{P}$ is
ergodic with respect to the {\em Wiener shift}    $\theta_t$:
\begin{equation}\label{NR2}
\theta_t\omega=\omega(\cdot+t)-\omega(t),\qquad\mbox{for
}\omega\in C_0({\Bbb R},U).
\end{equation}

A major example   of a random dynamical system is a random
differential equation. For example, let us  consider the following
evolution equation
 in some Hilbert space
\begin{equation}\label{eq-11}
\frac{du}{d\,t}=f(u,\theta_t\omega),\quad u(0)=x,
\end{equation}
over  some   metric dynamical system $(\Omega,{\mathcal F},{\Bbb
P},\theta)$. If (\ref{eq-11}) is well-posed for every $\om\in\Om$
and solutions $u(t,\om;x)$ depend measurably on $(t,\omega,x)$,
then the operator
\[
\phi\; :\; (t,\omega,x)\to u(t,\omega; x)
\]
 defines a random dynamical system (cocycle) $\phi$.
For  detailed presentation of random dynamical systems we refer to
the monograph by  Arnold \cite{Arn98}.
\medskip\par\noindent

Motivated by deterministic dynamical systems we introduce several
useful concepts from the theory of random dynamical systems.

A closed set $B(\omega)$,  depending on $\omega$,  in a separable
Hilbert space $H$
 is called random
if the distance mapping $\omega\to\sup_{x\in B(\omega)}\|x-y\|_H$ is a random variable for any $y\in H$.\\
A random dynamical system is called  {\em dissipative} if there
exists  a random set $B$   that is bounded for any $\omega$ and
that is absorbing: for  any random variable $x(\omega)\in H$ there
exists a $t_{x}(\omega)>0$ such that if $t\ge t_{x}(\omega)$,
then
\[
\phi(t,\omega,x(\omega))\in B(\theta_t\omega).
\]
 In the deterministic case ($\phi$ is independent of $\omega$)
the last relation coincides with the definition of an absorbing
set. In the case of partial differential equations of {\em
parabolic type}, due to the smoothing property, it is usually
possible to prove that a dissipative system possesses compact
invariant absorbing sets. For more details,  see  Temam
\cite{Tem97}, page 22f. Hence for a system of  stochastically
forced {\em parabolic} partial differential equations, such as
the stochastic two-layer fluid system introduced in the last
section, we usually consider the random set  $B(\omega)$ to be
compact. In addition, we will assume that $B(\omega)$ is  forward
invariant:
\[
\phi(t,\omega,B(\omega))\subset B(\theta_t\omega) , \; t>0.
\]
In the following we also need a concept of {\em tempered random
variables}. A random variable $x$ is called tempered if
\[
t\to |x(\theta_t\omega)|
\]
is   subexponentially growing:
\[
\limsup_{t\to\pm\infty}\frac{\log^+|x(\theta_t\omega)|}{|t|}=0
\quad\mbox{a.s.}
\]
where \[log^+(x) = max\{0,log(x)\}\]
This  technical  condition is
not a very strong restriction because the only alternative is that
the above $\limsup$ is $\infty$,
 which describes the degenerate case of stationarity;
see Arnold \cite{Arn98}, page 164 f.\\

\section{Cocycle property }

In this section we will show that the coupled fluid model
(\ref{q1}) defines a random dynamical system (cocycle property).
The cocycle property essentially comes from the well-posedness.\\
Now we   re-formulate the model such that appropriate tools of the
theory of random dynamical system can be applied to analyse the
coupled atmosphere-ocean model under a random wind forcing.
For the following we need some tools from the theory of partial differential equations.\\
Let $H^1(D)$ be the Sobolev space of functions on $D$ with first
generalized derivative in $L_2(D)$, the function space of square
integrable functions on $D$ with norm and  inner product
\[
\|u\|_{L_2}=\left(\int_D|u(x)|^2dD\right)^\frac{1}{2},\quad
(u,v)_{L_2}=\int_Du(x)v(x)dD,\quad u,\,v\in L_2(D).
\]
The space $H^1(D)$ is equipped with the norm
\[
\|u\|_{H^1}=\|u\|_{L_2}+\|\partial_y u\|_{L_2}+\|\partial_z
u\|_{L_2}.
\]
Motivated by the zero-boundary conditions of $q$ we also introduce
the space $\Hsob$ which contains roughly speaking  functions which
are zero on the boundary $\partial D$ of $D$. This space can be
equipped with the norm
\begin{equation}\label{eqno}
\|u\|_{\Hsob}=\|\partial_y u\|_{L_2}+\|\partial_z u\|_{L_2}.
\end{equation}

Another Sobolev space is given by $\dot{H}^1(D)$ which is a
subspace of $H^1(D)$ consisting of functions $u$ such that $\int_D
u dxdy=0$. A norm equivalent to the $H^1$-norm on $\dot{H}^1(D)$
is given by the right hand side of (\ref{eqno}). For functions in
$L_2(D)$
having this property we will write $\dot{L}_2(D)$.\\

For convenience, we introduce the vector notation for unknown
geophysical quantities  $u=(\Theta  ,q ,T)$

We now take the linear differential operator from (\ref{q1})

\[
A(u) = \left(
\begin{array}{l}
- \Delta\Theta + (1 + b(x,y))\Theta \\
- \nu \Delta q \\
- \Delta T\\
\end{array}
\right)
\]

Recall that the function $0<  b(y) < 1$. ${\sl A}$ is defined on
functions that are sufficiently smooth. We also have the following
boundary conditions from \S.2

\[
   \frac{\partial\Theta}{\partial n}  = \frac{\partial T }{\partial n} = 0,
   \]

\[
   \psi|{\partial D}=0,  \;  q|{\partial D} = 0.
   \]

We introduce the phase space for our geophysical quantities

\begin{align*}
\ H&=  L_2(D)\times L_2(D) \times L_2(D)\\
\ V&=  H^1(D)\times \Hsob \times \dot{H}^1(D).
\end{align*}

After this preparation, we are able to write our problem as a
stochastic evolution equation.
\\
Let $\dot{w}$ be a noise on $L_2(D)$ with finite energy  given by
the covariance operator  $Q$ of the Wiener process $w(t)$ which is
defined on a probability space $(\Omega,\mathcal{F},\mathbb{P})$.
For the vector
\[
W=(w,0,0),
\]
we rewrite the coupled atmosphere-ocean system (\ref{q1}) as a
stochastic evolution equation on   $V^\prime$:
\begin{equation}\label{eq19}
\frac{du}{dt}+Au=F(u)+\dot{W},\qquad u(0)=u_0\in H,
\end{equation}
where $\dot{w}$ is a white noise as the generalized temporal
derivative of a Wiener process  $w$  with continuous trajectories
on $\mathbb{R}$  and with values in $L_2(D)$. And F(u) is defined
as the follow:
$$
F(u) =\left(
\begin{array}{l}
-a + S_a(x,y) - b(y)S_o(x,y) + b(y)T(x,y) \\
- r q + PrRa\partial_y T - J(\psi, q + \beta  y )\\
- J(T, \psi)\\
\end{array}
\right)
$$

Sufficient for this regularity is that the trace of the covariance
is finite with respect to the space $L_2(D)$: ${\rm tr}_{L_2}Q <
\infty$. In particular, we can choose the canonical probability
space where the set of elementary events $\Omega$ consists of the
paths of $w$ and the probability measure $\mathbb{P}$ is the
Wiener measure with respect to covariance $Q$.
\\

In the following, we need  a stationary Ornstein-Uhlenbeck process
solving the linear  stochastic equation on $D$
\begin{equation}\label{eq20}
\frac{dz}{dt} + A_1z = \dot{w}
\end{equation}
where $A_1 = - \Delta + ( 1 + b(y))$ is the linear operator with
  zero Neumann boundary condition at $\partial D$ and with zero
initial condition.

\begin{lemma}
Suppose that the covariance $Q$ has a finite trace : ${\rm tr}_{L_2} Q<\infty$. Then
(\ref{eq20}) has a unique stationary solution generated by
\[
(t,\omega)\to z(\theta_t\omega).
\]
Moreover,  $Z(\omega)=(z(\omega),0,0)$ is a random variable in
$V$.
\end{lemma}

For the proof we refer to Da Prato and Zabczyk \cite{DaPZab92}, Chapter 5,
or Chueshov and Scheutzow \cite{ChuScheu01}.
\\

For our calculations it will be appropriate to transform (\ref{eq19})
into a differential equation without white noise but with random coefficients.
We  set
\begin{equation}\label{eq19a}
v:= u - Z
\end{equation}
Thus we obtain a random differential equation in $V^\prime$
\begin{equation}\label{eq18}
\frac{dv}{dt} + Av = F(v + Z(\theta_t\omega)),\qquad v(0) = v_0\in
H.
\end{equation}
Equivalently, we can formulate the  equation (\ref{eq18}) using test functions
\begin{equation*}
\frac{d}{dt}(v(t), \zeta) + a(v(t), \zeta) =
(F(v(t)+Z(\theta_t\omega)), \zeta)\quad \text{for all } \zeta\in
V.
\end{equation*}

We have obtained a differential equation without white noise but
with random coefficients.
Such a differential equation can be treated sample-wise
for {\em any} sample $\omega$. Hence it
  is simpler to consider (\ref{eq18})
  than to study the stochastic differential equation (\ref{eq19})
  directly.
We are looking for solutions in
\begin{equation*}
v\in C([0, \tau]; H)\cap L^2(0, \tau; V),
\end{equation*}
for all $\tau>0$. If we can solve this equation then $u:=v+Z$
defines a solution version of (\ref{eq19}). For the well posedness
of the problem we now have the following result.

\begin{theorem}\label{tEX}
({\bf Well-Posedness})
For any time $\tau>0$, there exists a unique solution of (\ref{eq18})
in $C([0,\tau];H)\cap L_2(0,\tau;V)$. In particular, the solution mapping
\[
{{\mathbb{R}}}^+\times \Omega\times H\ni(t,\omega,v_0)\to v(t)\in H
\]

is measurable in its arguments and the solution mapping
$H\ni v_0\to v(t)\in H$ is continuous.\\
\end{theorem}
\begin{proof}
By the properties of $A$ and $F$ the random   differential
equation  (\ref{eq18})  is   essentially similar  to the 2
dimensional  Navier Stokes equation.  Hence we have existence and
uniqueness and the above regularity assertions.
\end{proof}

On account of the transformation (\ref{eq19a}),  we find that (\ref{eq19}) also has a unique solution.

Since the solution mapping
\[
{\mathbb{R}}^+ \times\Omega\times H\ni (t, \omega, v_0)\to v(t,
\omega, v_0) =: \phi(t, \omega, v_0)\in H
\]
is well defined,  we can introduce a random dynamical system. On $\Omega$ we can define
a shift operator $\theta_t$ on the paths of the Wiener process
that pushes our noise:
\[
w(\cdot, \theta_t\omega) = w(\cdot + t,\omega) - w(t, \omega)\quad
\text{for }t\in{\mathbb{R}}
\]
which is called the {\em Wiener shift}.
Then $\{\theta_t\}_{t\in{\mathbb{R}}}$ forms a flow which is ergodic for the probability measure
${\mathbb{P}}$.
The properties of the solution mapping cause  the following relations
\begin{align*}
&\phi(t+\tau,\omega,u)=\phi(t,\theta_\tau\omega,\phi(\tau,\omega,u))\quad\text{for }t,\,\tau\ge 0\\
&\phi(0,\omega,u)=u
\end{align*}
for any $\omega\in\Omega$ and $u\in H$. This property is called the cocycle property of $\phi$
which is important to study the dynamics of random systems. It is a generalization
of the semigroup property. The cocycle $\phi$ together with the flow $\theta$ forms
a {\em random dynamical system}.

\section{Dissipativity}

In this section we are going to show that   the coupled
atmosphere-ocean system (\ref{q1})
 is dissipative,  in the sense that it has an absorbing  (random) set.
 This definition has been used for deterministic systems
\cite{Tem97}.
This means that the solution vector
$v$ is contained in a particular region of the phase space $H$ after a sufficiently long time.
Dissipativity will be very important for understanding  the
asymptotic dynamics of the system. This dissipativity will give us
estimate of the atmospheric temperature evolution under oceanic feedback.
Dynamical properties that follow  from this dissipativity will be considered in the next section. In particular, we will show that  the coupled  atmosphere-ocean system
has a random attractor, has finite degree of freedom,  and is ergodic
under suitable conditions.\\

Let
\begin{equation} \label{eqE} \tilde v=(\tilde \Theta,T),\quad
\tilde\Theta=\Theta-z.
\end{equation}
 By integration by parts and Poincare inequality, we have
\begin{align}\label{eqE1}
\begin{split}
\frac{d}{dt}(\frac{\lambda_0}{2}\|\tilde\Theta\|_{L^2}^2 +
\|T\|_{L^2}^2) \le & - 2\|\nabla\Theta\|_{L_2}^2 - \|\nabla
T\|_{L^2}^2\\
 - \frac12\|\Theta\|_{L^2}^2 &-
\frac{\lambda_0}{2}\|T\|_{L^2}^2 + 6a^2 + 6\|S_a\|_{L^2}^2 +
6\|S_o\|_{L^2}^2.
\end{split}
\end{align}

Here Poincare inequality $\lambda_0\|u\|^2 \le \|\nabla
u\|^2$$(\lambda_0 > 0)$ is used. Here and in the following we
stress  that $0 < b(y) < 1$.

For $q$, we have
\begin{equation}
\frac{d}{2dt}\|q\|_{L^2}^2 = - \nu\|\nabla q\|_{L^2}^2 -
r\|q\|_{L^2}^2 + PrRa\int_D\partial_y Tq dxdy - \beta\int_D
\frac{\partial\psi}{\partial x}q dxdy.
\end{equation}
Under the following condition
 \begin{equation}  4\nu r > \frac{\beta^2|l|^2}{\pi^2},
\label{condition}
   \end{equation}
 followed the argument of
\cite{gaodua03}, we have for some positive constant $\alpha > 0$
$$\frac{d}{2dt}\|q\|_{L^2}^2 \le -
\alpha\|q\|_{L^2}^2  + PrRa\int_D\partial_y Tq dxdy,$$ by the
Cauchy-Schwarz inequality, we get
\begin{equation}\frac{d}{dt}\|q\|_{L^2}^2 \le -
\alpha\|q\|_{L^2}^2  + \frac{Pr^2Ra^2}{\alpha}\|\nabla
T\|_{L^2}^2.
\end{equation}

Collecting all these estimates, we have  for some positive
constants $\alpha_1$ and $C_1$
$$
\frac{d}{dt}(\frac{\lambda_0}{2}\|\tilde\Theta\|_{L^2}^2 +
\|T\|_{L^2}^2+ \frac{\alpha\lambda_0}{2Pr^2Ra^2}\|q\|_{L^2}^2) +
\alpha_2\|\nabla v\|^2
$$
$$\le -\alpha_1(\frac{\lambda_0}{2}\|\tilde\Theta\|_{L^2}^2 +
\|T\|_{L^2}^2+ \frac{\alpha\lambda_0}{2Pr^2Ra^2}\|q\|_{L^2}^2) +
C_1
$$


By the Gronwall inequality, we finally conclude that
\begin{equation}\label{estimate-v}\|v\|_{H}^2 \le C_2\|v(0)\|_{H}^2e^{-\alpha_1 t}
+ C_3,
\end{equation}
and
\begin{equation}\label{estimate-v1}\int_s^t \|\nabla v(\tau)\|_{H}^2d\tau \le C_4\|v(0)\|_{H}^2e^{-\alpha_1 s}
+ C_5(t - s), \;\mbox{for any}\;0 \le s \le t.
\end{equation}

We now get the dissipativity of $v$. Roughly speaking
dissipativity means that all trajectories of the system move to a
bounded set in the phase space. For a random system we have the
following version of dissipativity.
\begin{definition}\label{defA}
A random set $B=\{B(\omega)\}_{\omega\in \Omega}$ consisting of closed bounded sets $B(\omega)$
is called absorbing for a random dynamical system $\phi$ if we have for
any  random set $D=\{D(\omega)\}_{\omega\in\Omega},\,D(\omega)\in H$ bounded,
 such that $t\to \sup_{y\in D(\theta_t\omega)}\|y\|_H$ has a subexponential
growth for $t\to\pm\infty$
\begin{align}\label{eqA0}
\begin{split}
&\phi(t,\omega,D(\omega))\subset B(\theta_t\omega) \quad\text{for }t\ge t_0(D,\omega)\\
&\phi(t,\theta_{-t}\omega,D(\theta_{-t}\omega))\subset B(\omega)\quad\text{for }t\ge t_0(D,\omega).
\end{split}
\end{align}
$B$ is called forward invariant if
\[
\phi(t,\omega,u_0)\in B(\theta_t\omega)\quad \text{if }  u_0\in B(\omega)\quad \text{for }t\ge 0.
\]
\end{definition}

From (\ref{estimate-v}), we could get the existence of absorbing
set $B(\omega) = \{v\in H, \|v\|_H^2 \le 2C_3\}$. Suppose that
$t\to\|v_0(\theta_{-t}\omega)\|_{H}^2$ growths not faster than
subexponential. Then we have the subexponential growth for
$B(\omega)$, and we get
%

\begin{lemma}\label{lEX}
The random set $B(\omega) = B(0, R(\omega)) = \{v\in H, \|v\|_H^2
\le R = 2C_3\}$ is an absorbing and forward invariant set for the
random dynamical system generated by (\ref{eq18}).
\end{lemma}

For the applications in the next section we need that the elements which are contained in the
absorbing set satisfy a particular regularity. To this end we introduce the function space
\[
{\mathcal{H}}^s:=\{u\in H: \|u\|_s^2:=\|A^\frac{s}{2}u\|_H^2<\infty\}
\]
where $s\in {\mathbb{R}}$. The operator $A^s$ is the $s$-th power of the positive
and symmetric operator $A$.
Note that these spaces are embedded in the Slobodeckij spaces
$H^s,\, s>0$. The norm of these spaces is denoted by $\|\cdot\|_{H^s}$.
This norm can be found in Egorov and Shubin \cite{EgoShu91}, Page 118.
But we do not need this norm explicitly. We only mention that on ${\mathcal{H}}^s$ the norm
$\|\cdot\|_s$
of
$H^s$ is equivalent to the norm of ${\mathcal{H}}^s$ for $0 < s$, see \cite{LioMag68}.\\

Our goal is it to show that  $v(1,\omega,D)$ is  a bounded set in
${\mathcal{H}}^s$ for some $s>0$.
This property causes the complete continuity of the mapping $v(1,\omega,\cdot)$.
We now derive a differential inequality for $t\|v(t)\|_s^2$.
By the chain rule we have
\[
\frac{d}{dt}(t\|v(t)\|_s^2)=\|v(t)\|_s^2+t\frac{d}{dt}\|v(t)\|_s^2.
\]
Note that for the embedding  constant $c_{6,s}$ between
${\mathcal{H}}^s$ and $V$
\[
\int_0^t\|v\|_s^2ds\le c_{6,s}^2 \int_0^t\|v\|_V^2ds\quad
\text{for } s \le 1
\]
such that the left hand  side is bounded if the initial conditions $v_0$
are contained in a bounded set in $H$.
The second term in the above formula can be expressed as followed:
\[
t\frac{d}{dt}(A^\frac{s}{2}v,A^\frac{s}{2}v)_H =
2t(\frac{d}{dt}v,A^sv)_H = - 2t(Av,A^sv)_H +
2t(F(v+Z(\theta_t\omega)),A^sv)_H.
\]
We have
\[
(Av,A^sv)_H=\|A^{\frac{1}{2}+\frac{s}{2}}v\|_H=\|v\|_{1+s}^2.
\]
For terms which including Jacobian operator, we apply some
embedding theorems, see Temam \cite{Tem83} Page 12, then taking
$F_3$ as a example we have
 for a constant $c_{5}>0$
\begin{align*}
(F_3(v),\zeta)_H&\le
c_{6}\,\|v\|_{m_1+1}\|\psi\|_{m_2+1}\|\zeta\|_{m_3},\,\quad
\zeta\in H_{m_3}
\end{align*}
where $m_1 + m_2 + m_3 \ge 1$ and $0 \le m_i < 1$. Here we use
that $D$ is of dimension 2. We then have for $m_1 = 0,\,m_2 = s <
1$ and $m_3=1-s$
\[
|( F_3(v),A^sT)_H|\le c_{7}\|T\|_V\|\psi\|_{1+s}\|T\|_{1+s}.
\]
$\|\psi(t)\|_{1+s}$ is bounded by $c_1^\prime\|q(t)\|_H$ by the
definition of  $\psi(t)$ and $\|v(t)\|_{L_\infty(0,T;H)}<\infty$.
Hence we have for any $\eps>0$ a constant $c_{8}(\eps)$:
$$
(F_3(v(t)),A^s T(t))_H \le c_{8}(\eps)
\|q\|_{L_\infty(0,T;H)}^2\|T(t)\|_V^2 + \eps \|T(t)\|_{1+s}^2$$ $$
\le  c_{8}(\eps) \|v\|_{L_\infty(0,T;H)}^2\|v(t)\|_V^2 + \eps
\|v(t)\|_{1+s}^2,
$$
where $\eps$ is chosen sufficiently small. The other terms could
be estimated similarly, we omit it here.
Using (\ref{estimate-v1}), we could obtain $\|v(t,\omega,v_0)\|_s,
\quad 0< s < \frac14$ is bounded for $t_0 \le T < \infty, t_0 > 0$
if $v_0$ is contained in a bounded set. This allows us to write
down the main assertion with respect to the dissipativity of this
section.
\begin{theorem}\label{tA}
For the random dynamical system generated by (\ref{eq18}),
there exists a compact random  set $B=\{B(\omega)\}_{\omega\in \Omega}$
which satisfies Definition \ref{defA}.
\end{theorem}
We define
\begin{equation}\label{eqAB}
B(\omega)=\overline{\phi(1,\theta_{-1}\omega,B(0,R(\theta_{-1}\omega)))}\subset
{\mathcal{H}}^s,\quad 0< s < \frac14.
\end{equation}
In particular, ${\mathcal{H}}^s$ is compactly embedded in $H$.

\section{Random dynamics}

In this section we   analyse the dynamical behavior of  the
coupled atmosphere-ocean system (\ref{q1}). However, it will be
enough to analyse the transformed  random dynamical system
generated by (\ref{eq19}). By the transformation (\ref{eq19a}) we
can take over
all these qualitative properties to the system (\ref{eq19}).\\

We will consider following dynamical behavior: random attractors,
atmospheric temperature evolution under oceanic feedback, and
ergodicity.

We first consider random climatic attractors.
We recall the following  basic concept; see, for instance,
Flandoli and Schmalfu{\ss} \cite{FlaSchm95a}.

\begin{definition}
Let $\phi$ be a random dynamical
dynamical system.
A random set $A=\{A(\omega)\}_{\omega\in \Omega}$ consisting of
compact nonempty sets $A(\omega)$ is called random global attractor
if for any random bounded set $D$ we have for the limit in probability
\[
({\mathbb{P}})\lim_{t\to\infty} {\rm dist}_H(\phi(t,\omega,D(\omega)),A(\theta_t\omega)) = 0
\]
and
\[
\phi(t,\omega,A(\omega))=A(\theta_t\omega)
\]
any $t\ge 0$ and $\omega\in\Omega$.
\end{definition}

The essential long-time behavior of a random system
is captured by a random  attractor.
In the last section we showed that the dynamical system $\phi$ generated by
(\ref{eq19}) is dissipative which means that there exists a random set $B$ satisfying
(\ref{eqA0}). In addition,  this set is compact. We now recall and adapt the following
theorem  from \cite{FlaSchm95a}.

\begin{theorem}
Let $\phi$ be a random dynamical
dynamical system on the state space $H$ which is a separable Banach space such that
$x\to\phi(t,\omega,x)$ is continuous.
Suppose that $B$ is a set ensuring the dissipativity given in definition \ref{defA}.
In addition, $B$ has a subexponential growth
(see  Definition \ref{defA})  and is regular (compact).
Then the dynamical system   $\phi$ has
a random attractor.
\end{theorem}

This theorem can be applied to our random dynamical system
 $\phi$ generated by the stochastic differential equation (\ref{eq19}).
 Indeed, all the assumptions
are satisfied. The set $B$ is defined in Theorem \ref{tA}. Its subexponential growth
follows from $B(\omega)\subset B(0,R(\omega))$ where
the radius $R(\omega)$ has been introduced in the last section.
Note that $\phi$ is a {\em continuous} random dynamical system; see Theorem \ref{tEX}. Thus $\phi$ has a random attractor.
By the transformation (\ref{eq19a}), this is also true for the original
coupled atmosphere-ocean system.

\begin{coro}  \label{attractor}
({\bf Random Attractor})
 The coupled atmosphere-ocean system (\ref{q1})
has a random attractor.
\end{coro}

\bigskip

Now we consider random fixed point and ergodicity. We can do a
small modification of (\ref{q1}). This modification is given when
we replace $\Delta T$ by $\nu\Delta T$  where $\nu>0$ is
viscosity. Under particular assumptions about physical data in
(\ref{q1}) , we can show that the behavior of our dynamical system
is laminar. For a stochastic system,  this means that after a
relatively short time, all trajectories starting from different
initial states show almost the same dynamical behavior. This can
be seen easily if $\,a$,  $S_a$ and $S_o$ are zero, there is no
noise and $\nu$ is large.
We will show that a laminar behavior also appears when  $\,a$, $S_a$ and $S_o$  are small in some sense.\\
Mathematically speaking,  laminar behavior means that a random dynamical system has a unique   exponentially
attracting random fixed point.

\begin{definition}
A random variable $v^\ast:\Omega\to H$ is defined to be a random fixed point for a random dynamical system
if
\[
\phi(t,\omega,v^\ast(\omega))=v^\ast(\theta_t\omega)
\]
for $t\ge 0$ and $\omega\in\Omega$. A random fixed point
$v^\ast$ is called exponentially attracting
if
\[
\lim_{t\to\infty}\|\phi(t,\omega,x)-v^\ast(\theta_t\omega)\|_H=0
\]
for any $x\in H$ and $\omega\in \Omega$.
\end{definition}
Sufficient conditions for the existence of random fixed points
are given in Schmalfu{\ss}
\cite{Schm97a}. We here formulate a simpler version of this
theorem and it is appropriate  for our system here.

\begin{theorem} ({\bf Random Fixed Point Theorem})
Let $\phi$ be a random dynamical system and suppose that $B$ is a forward invariant complete set.
In addition, $B$ has a subexponential growth, see Definition \ref{defA}.
Suppose that the following contraction conditions holds:
\begin{equation}\label{eqKX}
\sup_{v_1\not= v_2\in B(\omega)}\frac{\|\phi(1,\omega,v_1)-\phi(1,\omega,v_2)\|_H}{\|v_1-v_2\|_H}\le k(\omega)
\end{equation}
where the expectation of $\log k$ denoted by ${\mathbb{E}}\log k<0$.
Then $\phi$ has a unique random fixed point in $B$ which is exponentially attracting.
\end{theorem}

This theorem can be considered as a random version of the Banach fixed point theorem.
The contraction condition is formulated in the mean for the right hand side of  (\ref{eqKX}).\\

\begin{theorem} \label{fpt}
Assume  that the physical data $|a|$, $\|S_a\|_{L^2}$,
$\|S_o\|_{L^2}$ and the trace of the covariance for the noise
${\rm tr}_HQ$ are sufficiently small, and that the viscosity $\nu$
is sufficiently large. Then the random dynamical system generated
by (\ref{eq19}) has a unique random fixed point in $B$.
\end{theorem}
Here we only give a short sketch of the proof.
Let us suppose for a while that $B$ is given by the ball $B(0,R)$ introduced in Lemma \ref{lEX}.
Suppose that the data in the assumption of the lemma are small and $\nu$ is large.
Then it follows that ${\mathbb{E}}R$ is also small.
To calculate the contraction condition we have to calculate
 $\|\phi(1,\omega,v_1(\omega))-\phi(1,\omega,v_2(\omega))\|_H^2$
for arbitrary random variables  $v_1,\,v_2\in B$. By the property
of Jacobian operator we have that
\[
\langle J(q_1,\psi_1)-J(q_1,\psi_1),q_1-q_2\rangle\le
c_{9}\|q_1-q_2\|_{W_2^1}^2+c_{10}\|q_1\|_{\Hsob}^2\|q_1-q_2\|_{L_2}^2
\]
where the constant $c_{10}$ can be chosen sufficiently small if
$\nu$ is large.
By direct estimate, we get
\begin{align*}
\frac{d}{dt}&\|\phi(t,\omega,v_1(\omega))-\phi(t,\omega,v_2(\omega))\|_H^2\\
&\le (- \alpha^\prime +
c_{10}\|\phi(t,\omega,v_2(\omega))\|_V^2)\|\phi(t,\omega,v_1(\omega))
- \phi(t,\omega,v_2(\omega))\|_H^2
\end{align*}
for some positive $\alpha^\prime $ depending on $\nu$.
From this inequality and the Gronwall lemma it follows that the contraction condition (\ref{eqKX}) is satisfied if
\[
{\mathbb{E}}\sup_{u_2\in
B(\omega)}c_{10}\int_0^1\|\phi(t,\omega,v_2)\|_V^2dt <
\alpha^\prime.
\]
But by the energy inequality this property is satisfied if the ${\mathbb{E}}R$ and
${\mathbb{E}}\|z\|_V^2$ is sufficiently small which follows from the assumptions.
\\

Let now $B$ be the random set defined in (\ref{eqAB}).
Since the set $B$ introduced in (\ref{eqAB}) is absorbing any state the fixed point
$v^\ast$ is contained in this $B$. In addition $v^\ast$ attracts {\em any} state
from $H$ and not only states from $B$.

\begin{coro}\label{corrf}
({\bf Unique Random Fixed Point}) Assume that the  physical data
$|a|$, $\|S_a\|_{L^2}$, $\|S_o\|_{L^2}$ and the trace of the
covariance for the noise ${\rm tr}_HQ$ are sufficiently small, and
that the viscosity $\nu$ is sufficiently large. Then, through the
transformation (\ref{eq19a}), the original  system (\ref{q1}) has
a unique exponentially attracting random fixed point
$u^\ast(\omega) = v^\ast(\omega) + Z(\omega)$, where $u = (\Theta,
q, T)$.
\end{coro}

The uniqueness of  this random fixed point implies {\em ergodicity}.
We will comment on this issue at the end of this section.

\bigskip

By the well-posedness  Theorem \ref{tEX}, we know that the
stochastic differential equation (\ref{eq19}) for the coupled
atmosphere-ocean system has  a unique solution.  The solution is a
Markov process.  We can define the associated Markov operators
$\mathcal T(t)$ for $t \geq 0$, as discussed in \cite{Schm91a,
Schm89}. Moreover, $\{{\mathcal T}(t)\}_{t\geq 0}$ forms a
semigroup.

Let $M^2$ be the set of probability distributions $\mu$ with finite energy,
i.e.,
\[
\int_H\|u\|_H^2d\mu(u)<\infty .
\]
Then the distribution of the solution $u(t)$ (at time $t$) of the stochastic
differential equation   (\ref{eq19}) is given
by
\[
{\mathcal T}(t)\mu_0,\
\]
where the distribution $\mu_0$ of the initial data is contained in $M^2$.

We note that
the expectation of the solution $\|u(t)\|_H^2$ can be expressed in terms of this
distribution ${\mathcal T}(t)\mu_0$:
\[
{\mathbb{E}}\|u(t)\|_H^2=\int_H \|u\|_H^2d{\mathcal T}(t)\mu_0.
\]
We can derive the following energy inequality in the mean, using
our earlier estimates in (\ref{estimate-v}) and
(\ref{estimate-v1}):

\begin{theorem}\label{tMar}
The dynamical quantity  $u=(\Theta, q, T)$ of the coupled
atmospheric-ocean system(1) satisfy the estimate
\[
{\mathbb{E}}\|u(t)\|_H^2+\alpha{\mathbb{E}}\int_0^t\|u(\tau)\|_V^2d\tau
\le {\mathbb{E}}\|u_0\|_H^2+t\,c_{11}+t\,{\rm tr}_{L_2}Q,
\]
where the positive constants $c_{11}$ and $\alpha$ depend    on
physical data  $a$,   $\|S_a\|_{L^2}$, $\|S_o\|_{L^2}$, $Pr$ and
$Ra$.
\end{theorem}

By the Gronwall inequality, we further obtain the following result
about the asymptotic mean-square estimate.

\begin{coro} \label{feedbacktheorem}
 ({\bf Asymptotic Mean-Square Estimate })
For the expectation of the dynamical quantity  $u=(\Theta, q, T)$
of the coupled atmospheric-ocean system (1), we have the
asymptotic estimate
 \[
\limsup_{t\to\infty}
{\mathbb{E}}\|u(t)\|_H^2=\limsup_{t\to\infty}\int_H\|u\|_H^2d{\mathcal
T}(t)\mu_0\le \frac{c_{11}+{\rm tr}_{L_2}Q}{c_{12}}
\]
if the initial distribution $\mu_0$ of the random initial
condition $u_0(\omega)$ is contained in $M^2$.  Here   $c_{12} >
0$ also depends on physical data. In particular, we have
asymptotic mean-square estimate for the atmospheric temperature
evolution under oceanic feedback
\begin{equation} \label{temperature}
 \limsup_{t\to\infty}
{\mathbb{E}}\| \Theta \|_H^2
 \le
\frac{c_{11}+{\rm tr}_{L_2}Q}{c_{12}}.
\end{equation}
Thus the atmospheric temperature   $\Theta(y,t)$, as modeled by
the coupled atmosphere-ocean system (\ref{q1}),
 is bounded  asymptotically   in mean-square norm in terms of
physical quantities such as
the trace of the covariance operator of the external noise,
  the earth's longwave radiative cooling coefficient
$a$, and  the empirical functions $\|S_a\|_{L^2}$ and
$\|S_o\|_{L^2}$ representing the latitudinal dependence of the
shortwave solar radiation, as well as the  Prandtl  number ${\rm
Pr}$ and  the Rayleigh number ${\rm Ra}$ for oceanic fluids.
\end{coro}

By the estimates of Theorem \ref{tMar},
 we are able to use the well known
Krylov-Bogolyubov procedure to conclude the existence
of invariant measures of the Markov semigroup.
\begin{coro}
The semigroup of Markov operators $\{{\mathcal T}(t)\}_{t\geq 0}$ possesses an
invariant distribution $\mu_i$
in $M^2$:
\[
{\mathcal T}(t)\mu_i=\mu_i\quad \text{for }t\ge 0.
\]
\end{coro}

In fact,  the limit points of
\[
\left\{
\frac1t\int_0^t{\mathcal  T}(\tau)\mu_0d\tau\right\}_{t\geq 0}
\]
for $t\to\infty$ are invariant distributions.
The existence of such limit points follows from the estimate
in Theorem \ref{tMar}.

In some situations, the invariant measure may be unique.
For example,
the  unique random fixed point  in  Corollary \ref{corrf} is
 defined by a random variable
$u^\ast(\omega) = v^\ast(\omega)+z(\omega)$. This random variable
corresponds to a unique invariant measure of the Markov semigroup.
More specifically,
 this unique  invariant measure  is the expectation of the Dirac measure
with the random variable
as the random mass point
\[
\mu_i={\mathbb{E}}\delta_{u^\ast(\omega)}.
\]
Because the uniqueness of   invariant measure implies   ergodicity
\cite{DaPZab96}, we conclude that the
  coupled atmosphere-ocean model
 (\ref{q1})
is ergodic under the   suitable conditions  in  Corollary \ref{corrf} for physical data and random noise.  We reformulate Corollary \ref{corrf}
as the following  ergodicity principle.

\begin{theorem}  \label{ergodic}
({\bf Ergodicity}) Assume that the physical data $|a|$,
$\|S_a\|_{L^2}$, $\|S_o\|_{L^2}$
 and the trace of the covariance for the noise ${\rm
tr}_HQ$ are sufficiently small, and that the viscosity $\nu$ is
sufficiently large. Then the coupled atmosphere-ocean system
(\ref{q1}) is ergodic, namely, for any observable of the coupled
atmosphere-ocean flows, its time average approximates the
statistical ensemble average, as long as the time interval is
sufficiently long.
\end{theorem}

\bigskip

{\bf Acknowledgement.} This work was partly supported by the NSF
Grants DMS-0209326 and DMS-0139073,   and a Grant  of the NNSF of
China. A part of this work was done while J. Duan was visiting
Moscow State University, Russia, during May 2004. H. Gao would
like to thank the Illinois Institute of Technology, Chicago,
 for the hospitality.

\end{document}